\documentclass[11pt]{amsart}
\usepackage{amsmath,a4wide,scalerel}
\usepackage{mathabx}
\usepackage{xcolor,graphicx}
\usepackage{mathrsfs}
\usepackage{pgfplots}
\usepackage{subcaption}
\usepackage[normalem]{ulem}
\usepackage{float}
\usepackage{ctable} 

\usepackage{hyperref}

\newcommand{\R}{{\mathbb R}}
\newcommand{\N}{{\mathbb N}}
\newcommand{\E}{{\mathbb E}}

\definecolor{dkgreen}{rgb}{0,0.6,0}

\newcommand{\LT}{{\rm LT}}

\usepackage{color}
\usepackage{dsfont}

\begin{document}

\title[]{Positivity-preserving schemes for some nonlinear stochastic PDEs}

\author{Charles-Edouard Br\'ehier}
			\address{Universit\'e de Pau et des Pays de l'Adour, E2S UPPA, CNRS, LMAP, Pau, France}
			\email{charles-edouard.brehier@univ-pau.fr}

\author{David Cohen}
              \address{Department of Mathematical Sciences,
              Chalmers University of Technology and University of Gothenburg, 41296~Gothenburg, Sweden}
              \email{\tt david.cohen@chalmers.se}

\author{Johan Ulander}
              \address{Department of Mathematical Sciences,
              Chalmers University of Technology and University of Gothenburg, 41296~Gothenburg, Sweden}
              \email{\tt johanul@chalmers.se}

\begin{abstract}
We introduce a positivity-preserving numerical scheme for a class of nonlinear stochastic heat equations driven by a
purely time-dependent Brownian motion.
The construction is inspired by a recent preprint by the authors where one-dimensional equations driven by space-time white noise are considered. The objective of this paper is to illustrate the properties of the proposed integrators in a different framework, by numerical experiments and by giving convergence results.
\end{abstract}

\maketitle
{\small\noindent

Stochastic heat equation. Splitting scheme. Positivity. Mean-square convergence.

\section{Introduction}

Designing and studying numerical methods for stochastic partial differential equations (SPDEs) is an active field of research since the middle of the 1990's, we refer to the monograph~\cite{MR3308418} and to the recent preprint~\cite{bcu23} for a review of the literature. Proving sharp strong and weak convergence rates is not the only matter of interest, it is also desirable to preserve qualitative properties of the solutions at the discrete level,  see the classical reference \cite{MR2840298}. In order to illustrate this aspect we consider the following class of nonlinear heat equations driven by a multiplicative one-dimensional standard Brownian motion (using a formal notation for the noise $\dot{\beta}(t)$)
\begin{equation}\label{eq:spdebm-intro}
\left\lbrace
\begin{aligned}
&\partial_t u(t,x) = \Delta u(t,x) + g(u(t,x))\dot{\beta}(t)~,\quad t>0,~x\in\mathcal{D}\\
&u(t,x)=0~,\quad t\ge 0,~x\in\partial\mathcal{D},\\
&u(0,x) = u_{0}(x)~,\quad x\in\mathcal{D},
\end{aligned}
\right.
\end{equation}
for $(t,x) \in [0,T] \times \mathcal{D}$, where $\mathcal{D}=(0,1)^d$, see Section~\ref{sec:setting} for details on the notation and a rigorous formulation, see equation~\eqref{eq:spdebm}.

The noise in~\eqref{eq:spdebm-intro} is purely time-dependent and is interpreted in the It\^o sense. The nonlinearity $g \colon \mathbb{R} \to \mathbb{R}$ is of class $\mathcal{C}^1$ with bounded derivative, and is assumed to satisfy the condition $g(0)=0$. The above SPDE has the following qualitative property which follows from a comparison principle argument (see also~\cite{MR3069916}): if the initial condition $u_{0} \geq 0$ is continuous and nonnegative on $[0,1]^d$, then almost surely one has $u(t,x)\ge 0$ for all $t\ge 0$ and $x\in[0,1]^d$. Such property has also been proved for instance in~\cite{MR1149348,MR1271224,MR3262487} for SPDEs driven by space-time white noise, and we refer to the preprint~\cite{bcu23} for further references.

While classical time integrators, such as the Euler--Maruyama scheme, the semi-implicit Euler--Maruyama scheme,
and the stochastic exponential Euler integrator do converge when applied to the SPDE~\eqref{eq:spdebm-intro},
they do not satisfy the positivity-property of the exact solution to the SPDE (see below for a numerical illustration). In order to fix this issue, we propose a positivity-preserving explicit scheme~\eqref{eq:scheme}, based on a Lie--Trotter splitting strategy. The main idea of a splitting strategy is to decompose the vector field of the problem in such
a way that the obtained subsystems are exactly (or easily) integrated, see the monographs \cite{MR2840298,MR3642447}. In this work, we only deal with the temporal discretization. A fully-discrete scheme is easily obtained by combining the proposed time integrator with a standard finite difference method, which also preserves the positivity of the solution, see Section~\ref{sec:num}. Let us mention the recent works~\cite{MR4449543} for the construction and analysis of positivity-preserving schemes for linear SPDEs driven by a finite number of Brownian motions. We refer to the preprint~\cite{bcu23} for references on positivity-preserving schemes for stochastic differential equations.

In this short paper, we first briefly provide the necessary background to study~\eqref{eq:spdebm-intro} (Section~\ref{sec:setting}). The construction and the properties of the proposed scheme are given in Section~\ref{sec:scheme}. We state without proof the following main results: the Lie--Trotter splitting scheme is positivity-preserving and it converges in the mean-square sense to the solution of~\eqref{eq:spdebm-intro} with strong rate of convergence $1/2$. In future works, it may be interesting to identify the weak rate of convergence for the proposed scheme. Finally, Section~\ref{sec:num} presents numerical experiments in order to illustrate the superiority of the proposed integrator compared with classical ones.

The construction of the proposed positivity-preserving scheme~\eqref{eq:scheme} follows the same strategy as in the recent preprint~\cite{bcu23} written by the authors, where the case of one-dimensional nonlinear stochastic heat equations driven by space-time white noise interpreted in the It\^o sense (using a formal notation for the noise $\dot{W}(t,x)$)
\begin{equation}\label{eq:spde-intro}
\left\lbrace
\begin{aligned}
&\partial_t u(t,x) = \partial_{xx}^2 u(t,x) + g(u(t,x))\dot{W}(t,x),\\
&u(t,0)=u(t,1)=0,\\
&u(0,x) = u_{0}(x),
\end{aligned}
\right.
\end{equation}
for $(t,x)\in [0,T]\times (0,1)$ is considered. Let us briefly compare the results of this short paper with those of~\cite{bcu23}. First, note that~\eqref{eq:spde-intro} needs to be considered on a one-dimensional domain since it is driven by space-time white noise, whereas~\eqref{eq:spdebm-intro} can be considered in arbitrary dimension. Another major difference is the regularity of solutions: the solutions of~\eqref{eq:spde-intro} are H\"older continuous with exponent $1/4-$ in time and $1/2-$ in space, whereas the solutions of~\eqref{eq:spdebm-intro} are H\"older continuous with exponent $1/2-$ in time and $1-$ in space. As a result, the order of convergence of the splitting scheme differs, this is why one obtains strong order of convergence $1/2$ in this paper. Finally, in~\cite{bcu23} it is necessary to deal with a fully-discrete scheme, and to impose CFL stability conditions to ensure boundedness of moments and convergence of the scheme, even if the linear part of the problem is solved exactly (by an exponential integrator). Both the analysis and the numerical experiments in~\cite{bcu23} show the importance of the CFL conditions. On the contrary, in this paper the time-step size can be freely chosen and we are even able to study the scheme in a semi-discrete framework. The numerical experiments in Section~\ref{sec:num} show that indeed CFL conditions are not needed for the discretization of~\eqref{eq:spdebm-intro} using the proposed time integrator.

\section{Setting}\label{sec:setting}

In this work, we consider the following nonlinear stochastic heat equation driven by a purely time-dependent Brownian motion, interpreted in the It\^o sense:
\begin{equation}\label{eq:spdebm}
\left\lbrace
\begin{aligned}
&\text d u(t,x) = \Delta u(t,x)\,\text dt + g(u(t,x))\,\text d\beta(t)~,\quad t>0 , x\in\mathcal{D},\\
&u(t,x)=0~,\quad t\ge 0,~x\in\partial\mathcal{D},\\
&u(0,x) = u_{0}(x)~,\quad x\in \mathcal{D},
\end{aligned}
\right.
\end{equation}
where the spatial domain is $\mathcal{D}=(0,1)^d$ and $\overline{\mathcal{D}}=[0,1]^d$, in arbitrary dimension $d\ge 1$. Above $\Delta=\partial_{x_1x_1}^2+\ldots+\partial_{x_dx_d}^2$ is the Laplace operator and homogeneous Dirichlet boundary conditions are imposed on $\partial\mathcal{D}$. The evolution is driven by a standard real-valued Brownian motion $\bigl(\beta(t)\bigr)_{t\ge 0}$ defined on a probability space $(\Omega,\mathcal{F},\mathbb{P})$ satisfying the usual conditions.

The initial value $u_0:\overline{\mathcal{D}}\to\R$ is assumed to be a bounded and Lipschitz continuous mapping, and to satisfy the homogeneous Dirichlet boundary conditions: $u_0(x)=0$ for all $x\in\partial\mathcal{D}$. The initial value is assumed to be deterministic.
For all $\alpha\in(0,1]$, introduce the norms
\[
\|v\|_0=\underset{x\in\overline{\mathcal{D}}}\sup~|v(x)|~,\quad \|v\|_\alpha=\|v\|_0+\underset{x_1,x_2\in\overline{\mathcal{D}}}\sup~\frac{|v(x_2)-v(x_1)|}{|x_2-x_1|^\alpha}
\]
for any $\alpha$-H\"older continuous mapping $v$.

The nonlinearity $g:\R\to\R$ is a mapping of class $\mathcal{C}^1$, and is assumed to have a bounded derivative and to satisfy the condition $g(0)=0$.

Under the conditions above, the stochastic partial differential equation~\eqref{eq:spdebm} admits a unique mild solution, given by the integral formulation
\begin{equation}\label{eq:mildspdebm}
u(t,x)=\int_{\mathcal{D}} G(t,x,y)u_0(y)\,\text dy +\int_{0}^{t}\int_{\mathcal{D}}G(t-s,x,y)g(u(s,y))\,\text dy\,\text d\beta(s),~t\ge 0,x\in\overline{\mathcal{D}},
\end{equation}
where $(t,x,y)\in(0,+\infty)\times\overline{\mathcal{D}}^2\mapsto G(t,x,y)$ denotes the fundamental solution of the heat equation with homogeneous boundary conditions on the domain $\mathcal{D}$. We refer for instance to~\cite{MR3222416,MR876085} for standard references on the analysis of stochastic partial differential equations.

As seen in the introduction, the exact solution $\bigl(u(t,x)\bigr)_{t\ge 0,x\in\overline{\mathcal{D}}}$ of the SPDE~\eqref{eq:spdebm} satisfies the following property: if $u_0(x)\ge 0$ for all $x\in\overline{\mathcal{D}}$, then almost surely, one has $u(t,x)\ge 0$ for all $(t,x)\in[0,T]\times\overline{\mathcal{D}}$. See~\cite{MR3069916} for a proof. For a sketch of an alternative proof using the consistent positivity-preserving scheme~\eqref{eq:scheme}, see the end of Section~\ref{sec:scheme}.

Note that it would be straightforward to generalize the results presented in this paper to SPDEs driven by noise of the type $g(t,x,u(t,x))\,\text d\beta(t)$, for sufficiently regular functions $g$ satisfying the condition $g(t,x,0)=0$ for all $t\ge 0,x\in\overline{\mathcal{D}}$. Furthermore, with appropriate minor modifications we could also consider SPDEs driven by a noise of the type $\sum_{k=1}^{K}g_k(u(t,x))\,\text d\beta_k(t)$, where $\beta_1,\ldots,\beta_K$ are independent standard real-valued Brownian motions and the functions $g_k$ are of class $\mathcal{C}^1$, have bounded first order derivatives and satisfy the condition $g_k(0)=0$.
One could also extend the analysis to systems of SPDEs, like in~\cite{bcu23}. In the sequel we only consider the SPDE~\eqref{eq:spdebm} for ease of presentation.

\section{Positivity-preserving integrator}\label{sec:scheme}

Let us now describe the proposed time integrator for the approximation of the solution of~\eqref{eq:spdebm}. Let $T\in(0,\infty)$ be given and define the time-step size $\tau=T/M$ where $M\in\N$ is an integer. Set $t_m=m\tau$ for all $m\in\{0,\ldots,M\}$ and define the increments of the Brownian motion $\delta\beta_m=\beta(t_{m+1})-\beta(t_m)$ for all $m\in\{0,\ldots,M-1\}$. Introduce the auxiliary bounded and continuous function $f:\R\to\R$ defined by
\[
f(v)=\frac{g(v)}{v}\mathds{1}_{v\neq 0}+g'(0)\mathds{1}_{v=0}.
\]

The numerical approximation $u_m^{\LT}(\cdot)$ of the solution $u(t_m,\cdot)$ at time $t_m$ is defined as follows: for all $m\in\{0,\ldots,M-1\}$ and $x\in\overline{\mathcal{D}}$,
\begin{equation}\label{eq:scheme}
u_{m+1}^{\LT}(x)=\int_{\mathcal{D}} G(\tau,x,y)\left(\exp\Bigl(f(u_{m}^{\LT}(y))\delta\beta_m-\frac{f(u_{m}^{\LT}(y))^2\tau}{2}\Bigr) \right)\,\text dy,
\end{equation}
with initial value $u_0^{\LT}=u_0$. The proposed scheme~\eqref{eq:scheme} is based on a Lie--Trotter splitting strategy: given $u_m^{\LT}$ for some $m\in\{0,\ldots,M-1\}$, the numerical solution $u_{m+1}^{\LT}$ is obtained by solving successively two subsystems on the time interval $[t_m,t_{m+1}]$:
\begin{itemize}
\item first, the family of linear It\^o stochastic differential equations
\begin{equation}\label{eq:scheme-subsystem1}
\text dv_{1,m}(t,x)=v_{1,m}(t,x)f(u_{m}^{\LT}(x))\,\text d\beta(t)~,\quad t\in[t_m,t_{m+1}],~x\in\overline{\mathcal{D}},
\end{equation}
with initial value $v_{1,m}(t_m,\cdot)=u_{m}^{\LT}(\cdot)$;
\item second, the linear deterministic partial differential equation
\begin{equation}\label{eq:scheme-subsystem2}
\left\lbrace
\begin{aligned}
&\text dv_{2,m}(t,x)=\Delta v_{2,m}(t,x)\,\text dt~,\quad t\in (t_m,t_{m+1}),~x\in\mathcal{D}\\
&v_{2,m}(t,0)=v_{2,m}(t,1)=0~,\quad t\in [t_m,t_{m+1}],~x\in\partial\mathcal{D},
\end{aligned}
\right.
\end{equation}
with initial value $v_{2,m}(t_m,\cdot)=v_{1,m}(t_{m+1},\cdot)$.
\end{itemize}
Indeed, the exact solutions of the subsystem~\eqref{eq:scheme-subsystem1} and~\eqref{eq:scheme-subsystem2} are given by the following expressions: for all $t\in[t_m,t_{m+1}]$ and $x\in\overline{\mathcal{D}}$, one has
\begin{equation}\label{eq:scheme-subsystems-solutions}
\begin{aligned}
&v_{1,m}(t,x)=\exp\left(f(u_{m}^{\LT}(x))\bigl(\beta(t)-\beta(t_m)\bigr)-\frac{f(u_{m}^{\LT}(x))^2 (t-t_m)}{2}\right)u_{m}^{\LT}(x),\\
&v_{2,m}(t,x)=\int_{\mathcal{D}}G(t-t_m,x,y)v_{2,m}(t_m,y)\,\text dy=\int_{\mathcal{D}}G(t-t_m,x,y)v_{1,m}(t_{m+1},y)\,\text dy
\end{aligned}
\end{equation}
and the numerical approximation is set to $u_{m+1}^{\LT}(x)=v_{2,m}(t_{m+1},x)$, as prescribed by the Lie--Trotter splitting strategy. Note that the scheme~\eqref{eq:scheme} is explicit.

It is worth mentioning that the proposed scheme~\eqref{eq:scheme} is exact when applied to the linear stochastic heat equation~\eqref{eq:spdebm} when $g(v)=v$: in that case $u_m^{\LT}=u(t_m,\cdot)$ for all $m\in\{0,\ldots,M\}$. This can easily be seen by a change of unknown: if $g(v)=v$, then $(t,x)\mapsto e^{-\beta(t)+\frac{t}{2}}u(t,x)$ is the solution of the deterministic linear heat equation. In the general case, the nonlinearity $f$ is frozen at the left-point of each subinterval $[t_m,t_{m+1}]$, which results in the linear SDEs~\eqref{eq:scheme-subsystem1} which can then be solved exactly using~\eqref{eq:scheme-subsystems-solutions}.

The main benefit of introducing the explicit splitting scheme~\eqref{eq:scheme} is the following property: if $u_0(x)\ge 0$ for all $x\in\overline{\mathcal{D}}$, then for any choice of the time-step size $\tau=T/M$, one has $u_m^{\LT}(x)\ge 0$, for all $m\in\{0,\ldots,M\}$, almost surely. This means that the scheme is positivity-preserving. Proving this property is straightforward: by the interpretation as a splitting scheme, it suffices to check that the two subsystems~\eqref{eq:scheme-subsystem1} and~\eqref{eq:scheme-subsystem2} are positivity-preserving. This is easily seen in the expressions~\eqref{eq:scheme-subsystems-solutions} of their solutions $v_{1,m}(t,x)$ and $v_{2,m}(t,x)$.

The positivity-preserving property of the scheme~\eqref{eq:scheme} is ensured by a careful discretization of the stochastic perturbation term of~\eqref{eq:spdebm}, and is not satisfied for standard integrators. For instance, the stochastic exponential Euler integrator
\begin{equation}
u_{m+1}^{\rm SEXP}(x)=\int_{\mathcal{D}}G(\tau,x,y)\Bigl(u_m^{\rm SEXP}(y)+g(u_m^{\rm SEXP}(y))\delta\beta_m\Bigr)\,\text dy
\end{equation}
is not positivity-preserving since the support of the Gaussian random variables $\delta\beta_m$ is the entire real line, see also the numerical experiments below.

Let us now state properties of the numerical scheme~\eqref{eq:scheme} in order to justify that it is consistent with the SPDE~\eqref{eq:spdebm} when the time-step size
$\tau$ tends to zero. Recall that the initial value $u_0$ is Lipschitz continuous and that $\|u_0\|_0$ and $\|u_0\|_\alpha$ are defined in Section~\ref{sec:setting}.

First, moment bounds are satisfied: for all $T\in(0,\infty)$, there exists $C_0(T)\in(0,+\infty)$ such that for any time-step size $\tau=T/M$, one has
\begin{equation}\label{eq:momentbounds}
\underset{0\le m\le M}\sup~\underset{x\in\overline{\mathcal{D}}}\sup~\E[|u_m^{\LT}(x)|^2]\le C_0(T)\|u_0\|_0^2.
\end{equation}
Second, one has the following strong convergence result: for all $T\in(0,\infty)$ and all $\alpha\in(0,1)$, there exists $C_\alpha(T)\in(0,+\infty)$ such that for any time-step size $\tau=T/M$, one has
\begin{equation}\label{eq:strongerror}
\underset{0\le m\le M}\sup~\underset{x\in\overline{\mathcal{D}}}\sup~\E[|u_m^{\LT}(x)-u(t_m,x)|^2]\le C_\alpha(T)\|u_0\|_\alpha^2~\tau^{\alpha}.
\end{equation}
The strong error estimate~\eqref{eq:strongerror} states that the proposed integrator converges in a mean-square sense with order $1/2$. This order of convergence is expected to be optimal in general, as will be illustrated by the numerical experiments below.

Providing detailed proofs of the moment bounds~\eqref{eq:momentbounds} and of the strong error estimate~\eqref{eq:strongerror} is out of the scope of this work. It is worth mentioning that combining the positivity-preserving property of the scheme~\eqref{eq:scheme} and the strong error estimate~\eqref{eq:strongerror} provides a proof of the positivity of the exact solutions of the SPDE~\eqref{eq:spdebm}.

\section{Numerical experiments}\label{sec:num}

In this section we numerically illustrate the properties of the proposed scheme~\eqref{eq:scheme} and compare it with existing methods. We put emphasis on preservation of positivity and on mean-square error estimates in order to exhibit the strong rate of convergence $1/2$ given in Section~\ref{sec:scheme} above.

The one-dimensional stochastic nonlinear heat equation~\eqref{eq:spdebm} is first discretized
in space by a centered finite difference approximation on a uniform grid, see for instance \cite{MR1644183}
(for problems driven by space-time white noise). Let $N\in\N$, define the mesh size $h=1/N$, and the grid points $x_n=nh$ for $0\le n\le N$. We use the convention that for any vector $v=\bigl(v_n\bigr)_{1\le n\le N-1}\in\R^{N-1}$, we append discrete homogeneous Dirichlet boundary conditions $v_{0}=0$ and $v_{N}=0$ when needed. The spatially discrete $\R^{N-1}$-valued stochastic process
$u^N(t)=\bigl(u_n^N(t)\bigr)_{1\le n\le N-1}$, for all $t\ge 0$, is thus defined as the solution to the $N-1$-dimensional stochastic differential equation
\begin{equation}\label{eq:spatialscheme}
\text du^N(t)=N^2D^Nu^N(t)\,\text dt+g(u^N(t))\,\text d\beta(t),
\end{equation}
with initial value $u^N(0)=\bigl(u_0^N\bigr)_{1\le n\le N-1}=\bigl(u_0(x_n)\bigr)_{1\le n\le N-1}$, and the $(N-1)\times(N-1)$ matrix $D^N$ is the standard matrix for the approximation of the Laplace operator with homogeneous Dirichlet boundary conditions. The solution $u^N(t)$ of~\eqref{eq:spatialscheme} is nonnegative for nonnegative initial value $u^N(0)$, since $-D^N$ satisfies a monotonicity property.

The system of stochastic differential equations~\eqref{eq:spatialscheme} is then discretized in time by the following integrators (we recall that $\tau=T/M$ denotes the time step size):
\begin{itemize}
\item the proposed Lie--Trotter splitting scheme~\eqref{eq:scheme} (denoted LT below)
\begin{equation}
u^{\LT}_{m+1}=e^{\tau N^2D^N}\left(\exp\Bigl(f(u_{m}^{\LT})\delta\beta_m-\frac{f(u_{m}^{\LT})^2\tau}{2}\Bigr) \right)
\end{equation}
\item the Euler--Maruyama scheme (denoted EM below)
\[
u^{\rm EM}_{m+1}=u^{\rm EM}_m+\tau N^2D^Nu^{\rm EM}_m+g(u^{\rm EM}_m)\Delta_m\beta,
\]
\item the semi-implicit Euler--Maruyama scheme (denoted SEM below)
\[
u^{\rm SEM}_{m+1}=u^{\rm SEM}_m+\tau N^2D^Nu^{\rm SEM}_{m+1}+g(u^{\rm SEM}_m)\Delta_m\beta,
\]
\item the stochastic exponential Euler integrator (denoted SEXP below)
\[
u^{\rm SEXP}_{m+1}=e^{\tau N^2D^N}\left(u^{\rm SEXP}_m+g(u^{\rm SEXP}_m)\Delta_m\beta\right).
\]
\end{itemize}

In the first numerical experiment, we illustrate the positivity-preserving property of the Lie--Trotter scheme (LT)
when applied to the time discretization of the stochastic heat equation~\eqref{eq:spdebm}
on the time interval $[0,2]$ with the following
multiplicative terms: $g(v)=\lambda v$, $g(v)=\lambda v/(1+v^2)$, $\lambda(\sin(v)+v)$, and $g(v)=\lambda\ln(1+v)$, where the real parameter $\lambda$ is introduced to modify the size of the noise.
We consider the following parameters:
$u_0(x)=\sin(\pi x)$, $\tau=2^{-5}$, $N=2^8$, $\lambda=2.5$ and compute $100$ realizations of each time integrators. The results
are presented in Table~\ref{tabBM}. The proposed scheme produces only nonnegative numerical solutions, which confirms the result stated in Section~\ref{sec:scheme}. On the contrary, the other integrators produce some solutions with negative values. This illustrates the superiority of the proposed scheme~\eqref{eq:scheme}.

\begin{table}[h]
\begin{center}
\begin{tabular}{|c c c c c|}
  \hline
  $g(v)$ & LT & EM & SEM & SEXP \\
  \specialrule{.1em}{.05em}{.05em}
  $2.5v$ & $100/100$ & $2/100$ & $47/100$ & $47/100$ \\
  \hline
  $2.5v/(1+v^2)$ & $100/100$ & $2/100$ & $49/100$ & $49/100$ \\
  \hline
  $2.5(\sin(v)+v)$ & $100/100$ & $0/100$ & $2/100$ & $2/100$ \\
  \hline
  $2.5\ln(1+v)$ & $100/100$ & $2/100$ & $50/100$ & $49/100$ \\
  \hline
\end{tabular}
\end{center}
\caption{Proportion of samples containing only positive values out of $100$ simulated sample paths for the time integrators: Lie--Trotter scheme (LT), Euler--Maruyama scheme (EM), semi-implicit Euler--Maruyama (SEM), and
stochastic exponential Euler scheme (SEXP). Time-step size: $\tau=2^{-5}$. Mesh size: $h=2^{-8}$.}
\label{tabBM}
\end{table}

In the second numerical experiment, we investigate the mean-square errors of the above time integrators in order to confirm the convergence result stated in Section~\ref{sec:scheme}.
We discretize the stochastic heat equation~\eqref{eq:spdebm} on the time interval $[0,0.5]$ with $g(v)=v$ and $g(v)=v/(1+v^2)$ and
initial value $u_0(x)=\sin(\pi x)$. The spatial discretization is again performed by a centered finite difference method with mesh size $h=2^{-8}$.
The temporal discretizations is done by the time integrators: LT, SEXP, and SEM. In this experiment the explicit EM integrator is not tested. Figure~\ref{fig:strongBM} presents, in a loglog plot, the mean-square errors
\[
\underset{0\le m\le M}\sup~\underset{0\le n\le N}\sup~\bigl(\E[|u_{m,n}^{\text{num}}-u^{\text{ref}}(t_m,x_n)|^2]\bigr)^{\frac12}
\]
measured for the time interval $[0,0.5]$. The time step sizes used for these experiments range from $\tau=2^{-4}$ to $\tau=2^{-16}$. The reference solution $u^{\text{ref}}$ is computed using the Lie--Trotter splitting scheme with $\tau=2^{-16}$. We use $150$ samples
to approximate the expectations. We have experimentally checked that the Monte Carlo error is negligible to observe mean-square convergence.
In the first plot of Figure~\ref{fig:strongBM}, one observes that if $g(v)=v$ then the Lie--Trotter splitting scheme produces the exact solution as explained in Section~\ref{sec:scheme}, while the other integrators have rate of convergence $1/2$. In the second plot of Figure~\ref{fig:strongBM}, one observes a rate of convergence $1/2$ in the mean-square error estimates for the three integrators. This confirms the convergence result stated in Section~\ref{sec:scheme}.

\begin{figure}[h]
\centering
\begin{subfigure}{.4\textwidth}
  \centering
  \includegraphics[width=\textwidth]{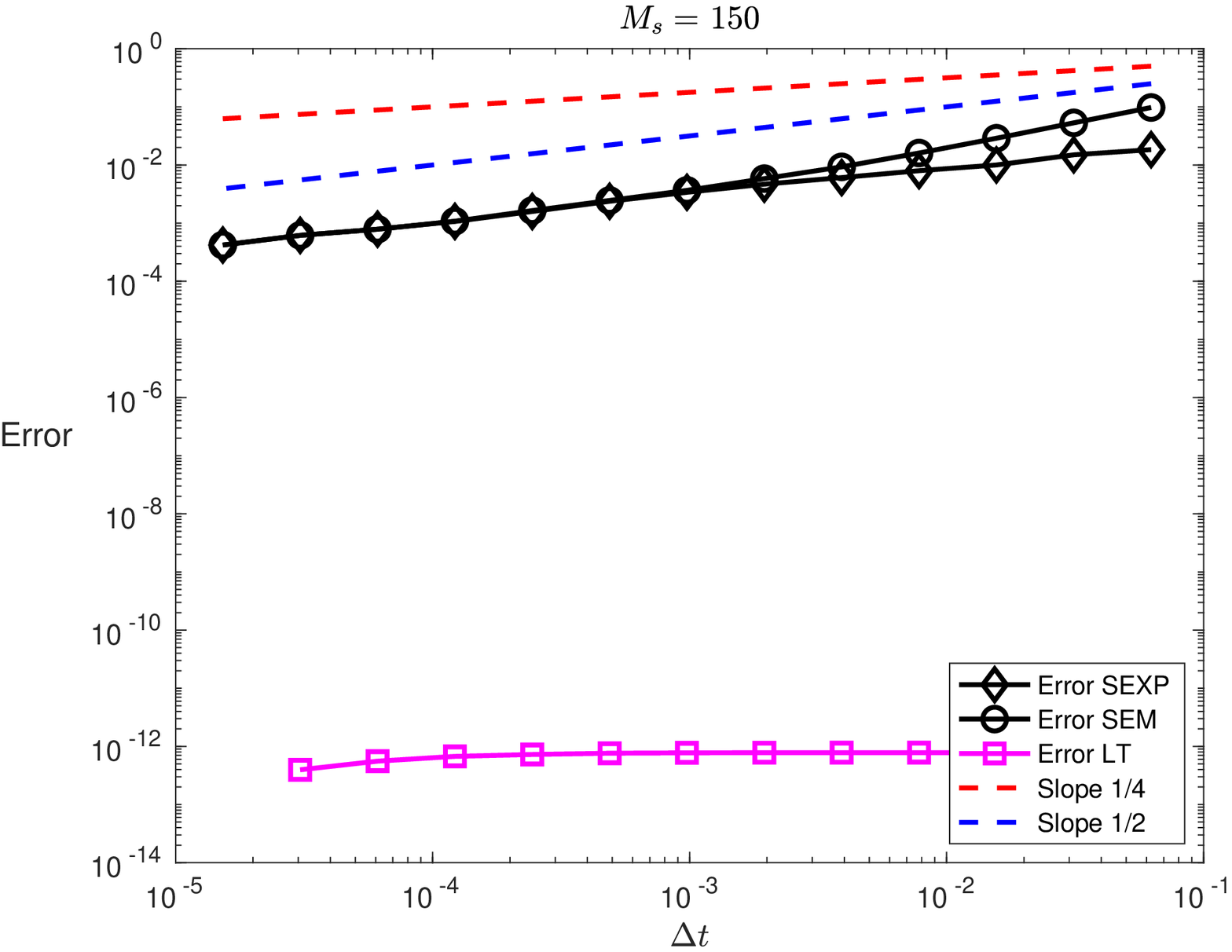}
  \caption{$g(v)=v$}
\end{subfigure}%
\begin{subfigure}{.4\textwidth}
  \centering
  \includegraphics[width=\textwidth]{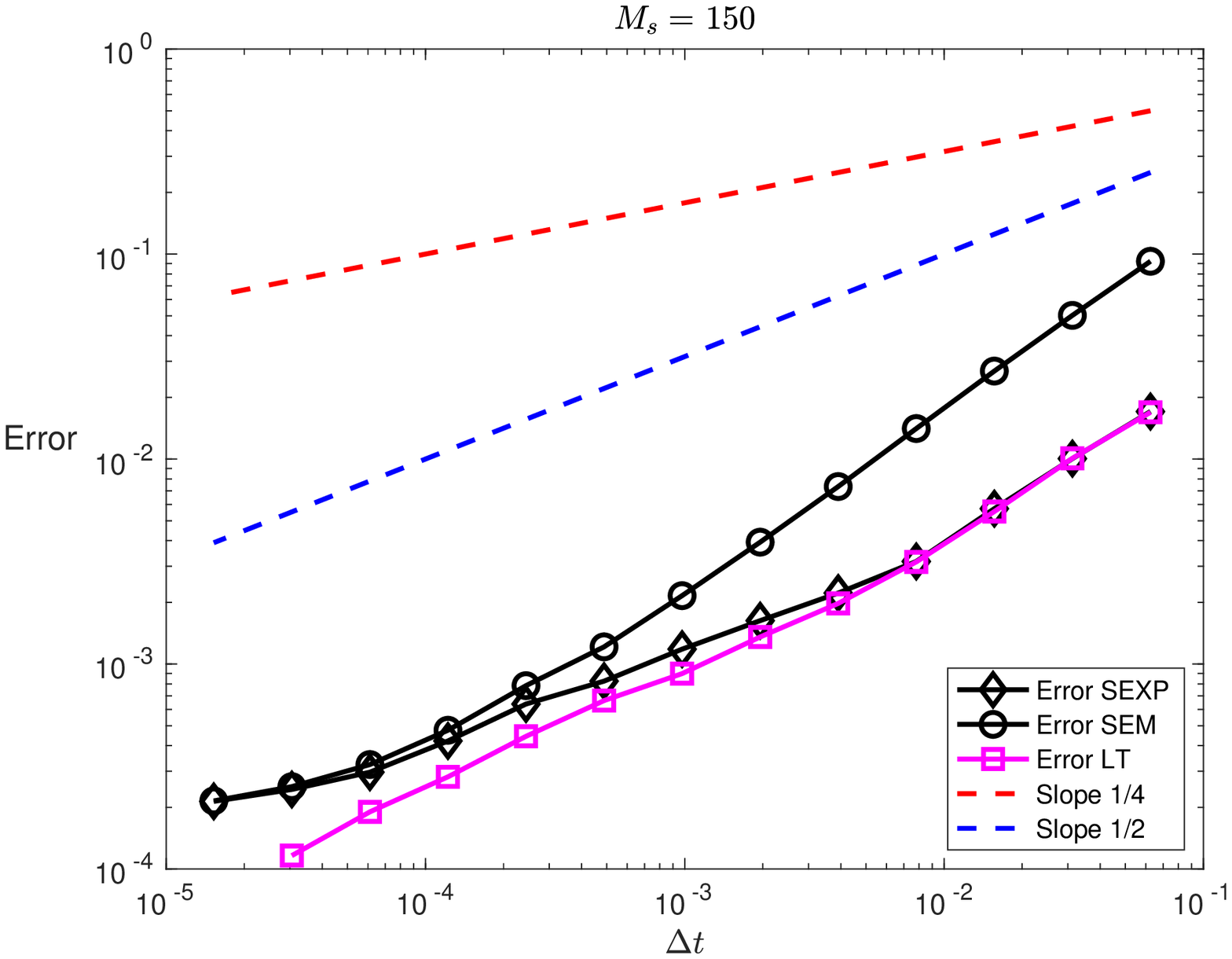}
  \caption{$g(v)=\frac{v}{(1+v^2)}$}
\end{subfigure}%

\caption{Mean-square errors of the splitting scheme (LT), the stochastic exponential Euler integrator (SEXP), and the semi-implicit Euler--Maruyama scheme (SEM). Mesh size $h=2^{-8}$ and average over $150$ samples.}
\label{fig:strongBM}
\end{figure}

Finally, we illustrate the fact that these error bounds are uniform in the spatial discretization. We compute the mean-square errors on the time interval $[0,0.5]$ of the Lie--Trotter splitting scheme when applied to the finite difference discretization of the stochastic heat equation with $g(v)=1.5v$, resp. $g(v)=1.5v/(1+v^2)$, and mesh sizes $h=2^{-4}, 2^{-6}, 2^{-8}, 2^{-10}$.
The time step sizes used for these experiments range from $\tau=2^{-4}$ to $\tau=2^{-16}$. The reference solutions are computed using the Lie--Trotter splitting scheme with $\tau=2^{-16}$. As above $150$ samples are used to approximate the expectations and the Monte Carlo error is negligible for the observation of the rates of convergence. These results are presented in Figure~\ref{fig:changeNgBM}. One observes that the error does not depend on the mesh size $h$. This is in sharp contrast to the observations from the preprint~\cite{bcu23} on the approximation of the equation~\eqref{eq:spde-intro} driven by space-time white noise, for which a CFL condition is required.

\begin{figure}[h]
\centering
\begin{subfigure}{.4\textwidth}
  \centering
  \includegraphics[width=\textwidth]{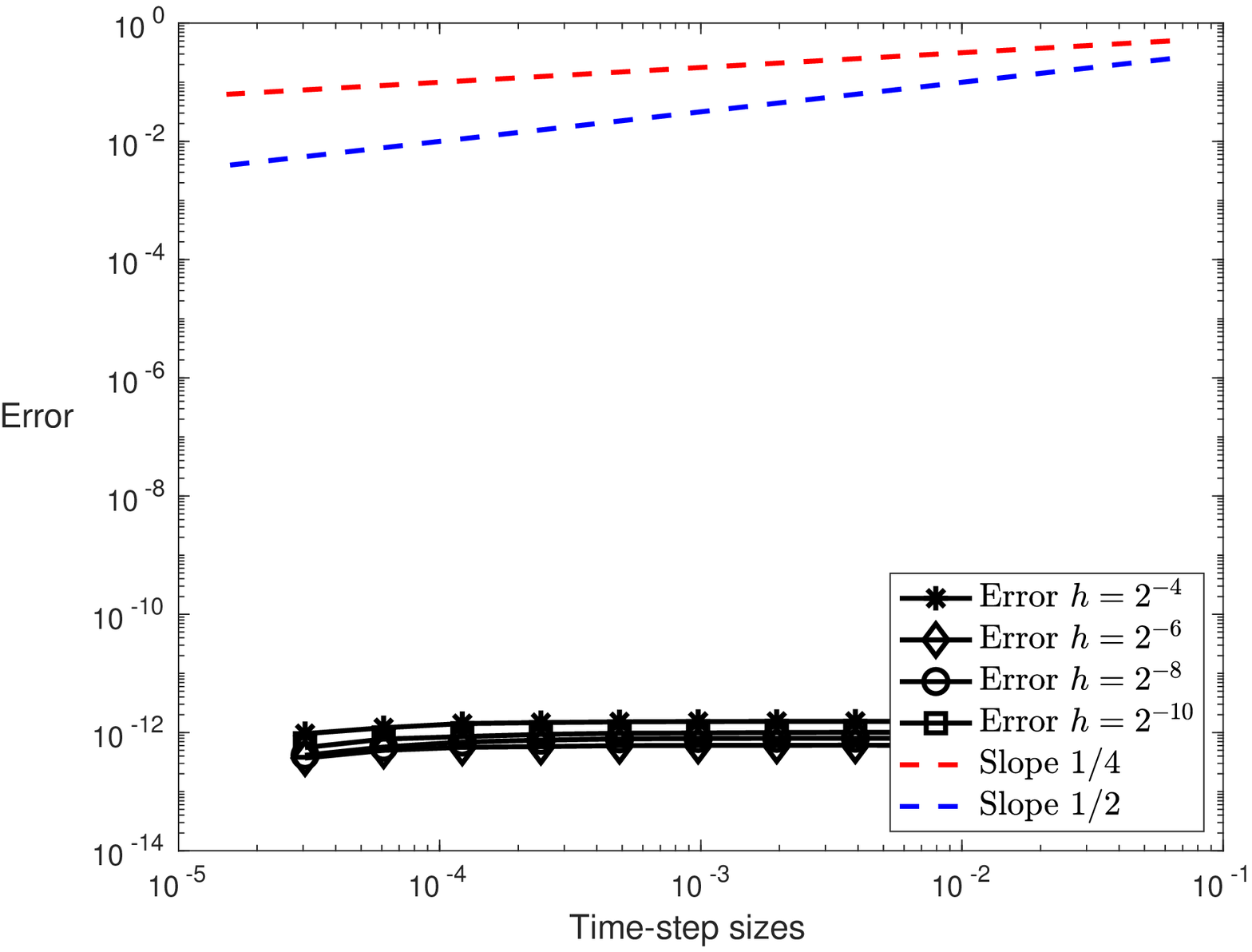}
  \caption{$g(v)=1.5v$}
\end{subfigure}%
\begin{subfigure}{.4\textwidth}
  \centering
  \includegraphics[width=\textwidth]{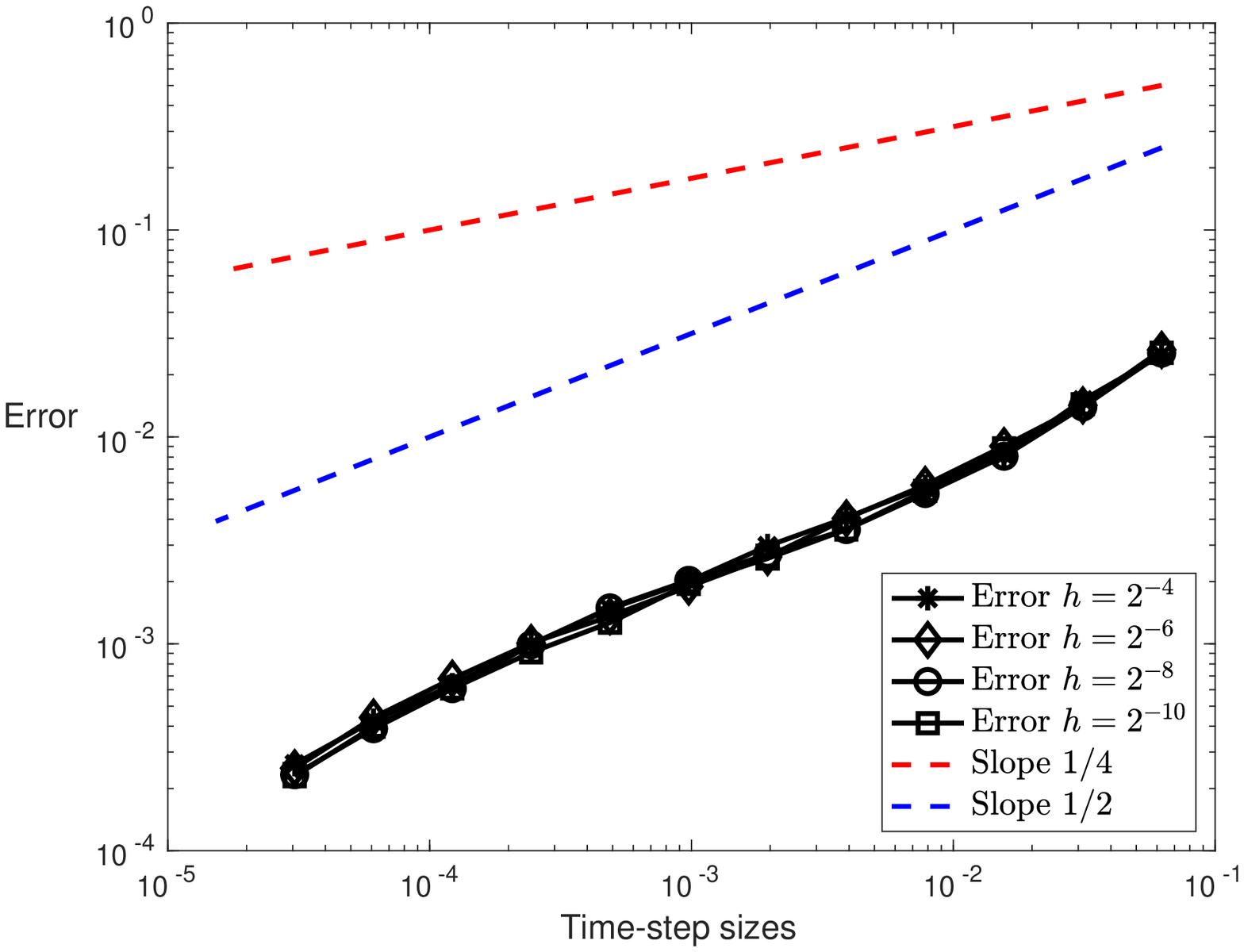}
  \caption{$g(v)=1.5\frac{v}{(1+v^2)}$}
\end{subfigure}%

\caption{Mean-square errors of the splitting scheme for several values of the spatial mesh size $h=2^{-4}, 2^{-6}, 2^{-8}, 2^{-10}$. Average over $150$ samples.}
\label{fig:changeNgBM}
\end{figure}

We conclude this paper with some numerical experiments in dimension $d=2$.

Let us first consider the stochastic heat equation~\eqref{eq:spdebm} on the time interval $[0,2]$
with initial value $u_0(x_1,x_2)=\sin(\pi x_1)\sin(\pi x_2)$ and with multiplicative terms:
$g(v)=2.5 v$, $g(v)=2.5 v/(1+v^2)$, $2.5(\sin(v)+v)$, and $g(v)=2.5\ln(1+v)$.
The discretization parameters are taken to be $\tau=2^{-5}$ and $h_{x_1}=h_{x_2}=2^{-4}$.
We compute $100$ realizations of each time integrators. The proportion of samples containing only positive values is presented in Table~\ref{tabBM2d}. One can again observe the superiority of the proposed Lie--Trotter splitting scheme.

\begin{table}[h]
\begin{center}
\begin{tabular}{|c c c c c|}
  \hline
  $g(v)$ & LT & EM & SEM & SEXP \\
  \specialrule{.1em}{.05em}{.05em}
  $2.5v$ & $100/100$ & $0/100$ & $47/100$ & $47/100$ \\
  \hline
  $2.5v/(1+v^2)$ & $100/100$ & $0/100$ & $48/100$ & $48/100$ \\
  \hline
  $2.5(\sin(v)+v)$ & $100/100$ & $0/100$ & $2/100$ & $2/100$ \\
  \hline
  $2.5\ln(1+v)$ & $100/100$ & $0/100$ & $46/100$ & $53/100$ \\
  \hline
\end{tabular}
\end{center}
\caption{SPDE in $2d$: Proportion of samples containing only positive values out of $100$ simulated sample paths for the time integrators: Lie--Trotter scheme (LT), Euler--Maruyama scheme (EM), semi-implicit Euler--Maruyama (SEM), and
stochastic exponential Euler scheme (SEXP). Time-step size: $\tau=2^{-5}$. Mesh sizes: $h_{x_1}=h_{x_2}=2^{-4}$.}
\label{tabBM2d}
\end{table}

Next, we compute the mean-square errors, measured for the time interval $[0,0.5]$, of the LT, SEXP and SEM integrators
when applied to the SPDE~\eqref{eq:spdebm} with $g(v)=v$ and $g(v)=v/(1+v^2)$ and
initial value $u_0(x_1,x_2)=\sin(\pi x_1)\sin(\pi x_2)$. The discretization parameters are:
$h_{x_1}=h_{x_2}=2^{-4}$ for the mesh sizes and the time-step size ranging from
$\tau=2^{-4}$ to $\tau=2^{-14}$. The reference solution $u_{\text{ref}}$ is computed using the Lie--Trotter splitting scheme with $\tau=2^{-14}$. We use $150$ samples to approximate the expectations.
The results are presented in Figure~\ref{fig:strongBM2d}. Again one observes that the Lie--Trotter splitting
scheme is exact for linear problems and has a rate of convergence $1/2$ in the mean-square sense.

\begin{figure}[h]
\centering
\begin{subfigure}{.4\textwidth}
  \centering
  \includegraphics[width=\textwidth]{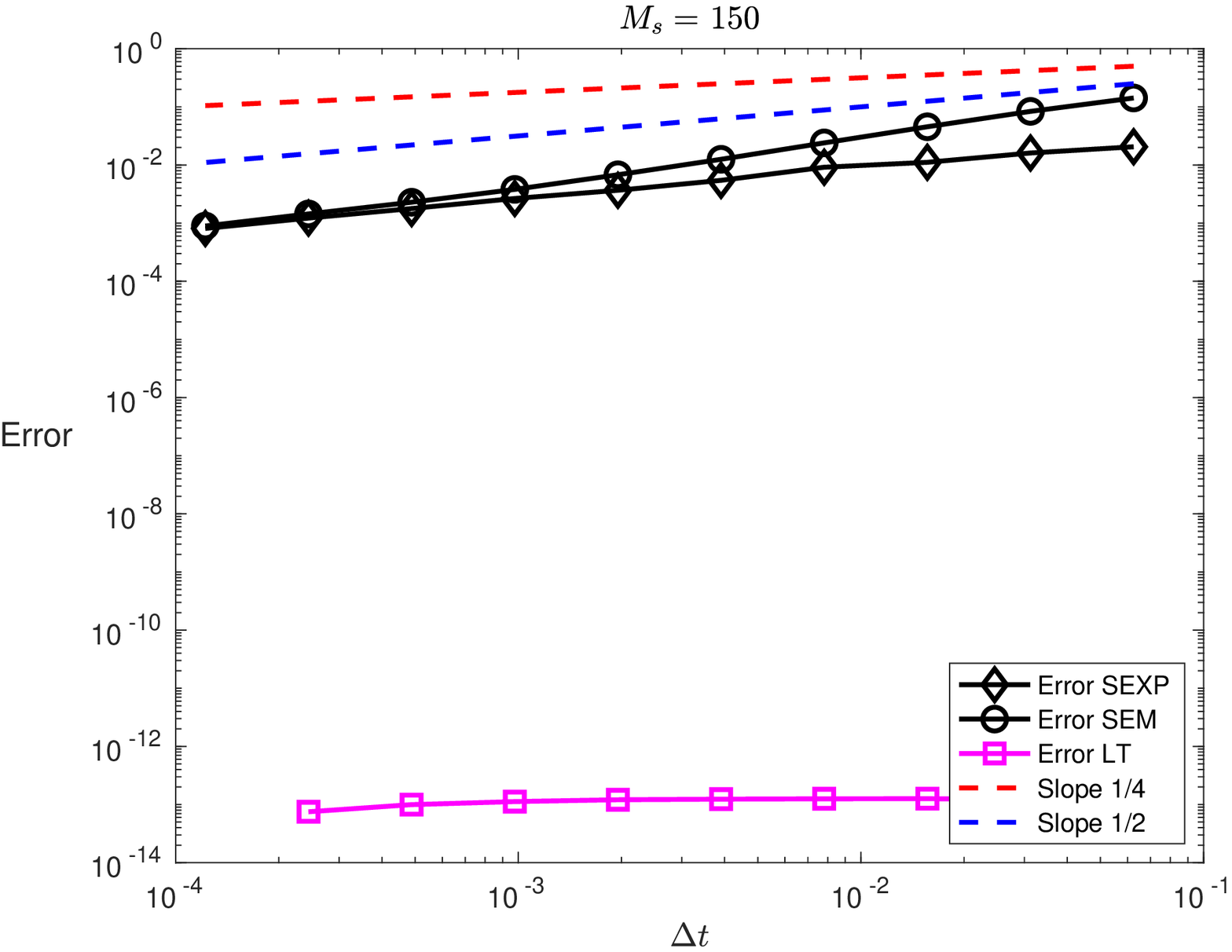}
  \caption{$g(v)=v$}
\end{subfigure}%
\begin{subfigure}{.4\textwidth}
  \centering
  \includegraphics[width=\textwidth]{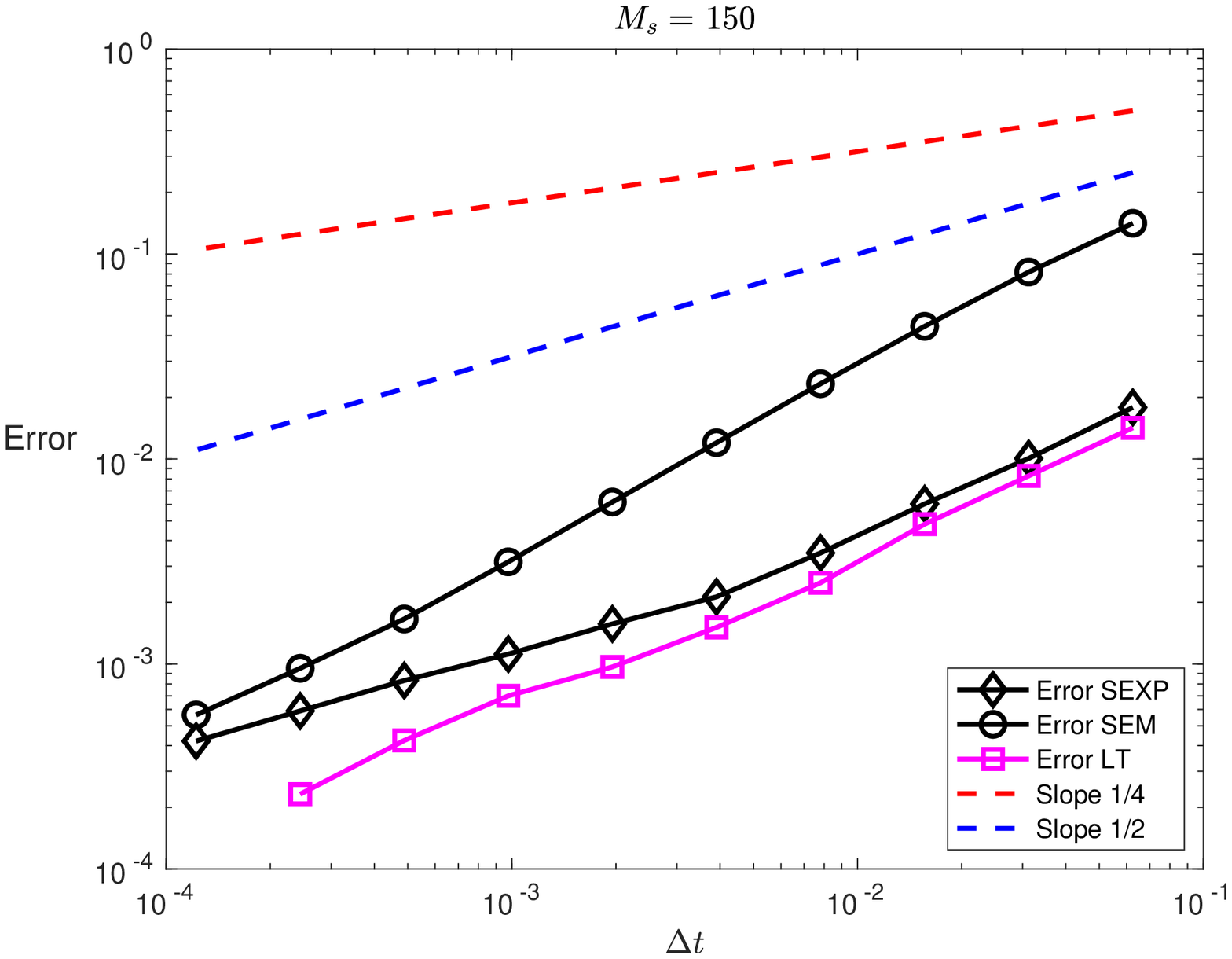}
  \caption{$g(v)=\frac{v}{(1+v^2)}$}
\end{subfigure}%

\caption{SPDE in $2d$: Mean-square errors of the splitting scheme (LT), the stochastic exponential Euler integrator (SEXP), and the semi-implicit Euler--Maruyama scheme (SEM). Mesh sizes $h_{x_1}=h_{x_2´}=2^{-4}$ and average over $150$ samples.}
\label{fig:strongBM2d}
\end{figure}

\section*{Acknowledgements}

We thank the referee for helpful comments on an earlier version of the paper. The work of CEB is partially supported by the project SIMALIN (ANR-19-CE40-0016) operated by the French National Research Agency.
The work of DC and JU is partially supported by the Swedish Research Council (VR) (projects nr. $2018-04443$).
The computations were performed on resources provided by the Swedish National Infrastructure
for Computing (SNIC) at HPC2N, Ume{\aa} University and at UPPMAX, Uppsala University.


\end{document}